\documentclass{article}

\usepackage{arxiv}

\usepackage[utf8]{inputenc} % allow utf-8 input
\usepackage[T1]{fontenc}    % use 8-bit T1 fonts
\usepackage{hyperref}       % hyperlinks
\usepackage{url}            % simple URL typesetting
\usepackage{booktabs}       % professional-quality tables
\usepackage{amsfonts}       % blackboard math symbols
\usepackage{nicefrac}       % compact symbols for 1/2, etc.
\usepackage{microtype}      % microtypography
\usepackage{lipsum}		% Can be removed after putting your text content
\usepackage{graphicx}
\usepackage{natbib}
\usepackage{doi}
\usepackage{placeins}

\usepackage{tikz}
\usepackage{amsmath}

\usepackage{mathtools,amsfonts,amssymb,bm}
\usepackage{amsthm}
\newtheorem{proposition}{Proposition}
\newtheorem{corollary}{Corollary}
\theoremstyle{remark}
\newtheorem*{remark}{Remark}
\usetikzlibrary{arrows.meta,positioning,calc}

\title{Switching Structure and Finite-Event Optima for One-Orbit Minimum-Fuel In-Plane HCW Maneuvers}

\author{ \href{https://orcid.org/0009-0009-2587-2474}{\includegraphics[scale=0.06]{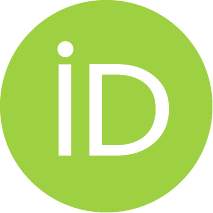}\hspace{1mm}Tarik Ziya Inci} \\
	Department of Astronautical Engineering\\
	Istanbul Technical University\\
	Istanbul, TURKEY \\
	\texttt{incit22@itu.edu.tr} \\
	\And
	\href{https://orcid.org/0000-0002-3268-1467}{\includegraphics[scale=0.06]{orcid.pdf}\hspace{1mm}Burak Boyacioglu} \\
    Center for Robotics and AI\\
	Middle East Technical University\\
	Ankara, TURKEY \\
	\texttt{bycgl@metu.edu.tr} \\
}

\renewcommand{\shorttitle}{Switching Structure of Minimum-Fuel HCW Maneuvers}

\hypersetup{
  pdftitle={Switching Structure and Finite-Event Optima for One-Orbit Minimum-Fuel In-Plane HCW Maneuvers},
  pdfsubject={math.OC, eess.SY},
  pdfauthor={Tarik Z. Inci, Burak Boyacioglu},
  pdfkeywords={Minimum-fuel control, Hill-Clohessy-Wiltshire dynamics,
    Spacecraft relative motion, Switching structure, Bang-off-bang control,
    Singular arcs},
}

\begin{document}
\maketitle

\begin{abstract}
	Under in-plane Hill--Clohessy--Wiltshire dynamics with independently bounded radial and along-track thrust, we prove that every target reachable within one orbit has a fixed-time minimum-fuel control with at most ten switches. We first show that regular optima have at most four radial and six along-track switches, coupled through a shared costate. Of the 35 candidate count pairs, 11 require minimum durations, one cannot be optimal, and 23 are attainable at arbitrarily short durations. We then convert singular-arc optima into bang/off optima with at most five switches by preserving terminal moments. We demonstrate terminal accuracy and whole-arc Pontryagin consistency through 55,296 numerical transfers.
\end{abstract}

% keywords can be removed

\section{Introduction}
\label{sec:introduction}

Spacecraft formation flying, rendezvous, inspection, servicing, and proximity operations require repeated planning of relative motion between prescribed local states \citep{DIMAURO2018311,LUO20141,ortolano2021autonomous,electronics12143034,OPROMOLLA2018342}. In these settings, the maneuvering spacecraft is often described with respect to a nearby reference vehicle or circular reference orbit, and the local-vertical local-horizontal (LVLH) frame provides a natural coordinate system for expressing the relative position and velocity. For circular reference orbits and  separations small relative to the orbital radius, the Hill--Clohessy--Wiltshire (HCW) equations give a standard linear model for relative motion with respect to the reference orbit or spacecraft \citep{cw1960,vallado2022}. Because the resulting state transition matrix is available in closed form, the HCW model remains attractive for local guidance and trajectory generation, especially when a higher-level planner repeatedly provides nearby target states to be reached over prescribed maneuver durations.

For such missions, propellant usage is a central design consideration because it directly affects mission lifetime and the number of maneuvers that can be performed \citep{YOO20131476}. This consideration is particularly important for small spacecraft and distributed spacecraft missions, where onboard resources and propulsion authority are often limited \citep{miller2021survey}. Consequently, local maneuver planning is commonly formulated with bounded control inputs and a fuel-related objective~\citep{carterbrient1992,gazzino2018,pontani2015}. In the present work, we consider a fixed-final-time transfer between prescribed relative states, with acceleration inputs bounded component-wise and fuel consumption represented by the integral of the absolute control effort. This formulation naturally leads to a minimum-$L_1$ control problem under HCW dynamics, where the objective is to reach the desired terminal state while minimizing the total thrusting effort over the maneuver interval.

Indirect optimal-control methods have been used extensively to study low-thrust relative-motion problems under the HCW model~\citep{force1,force2,carterbrient1991,pontani2015}. A useful feature of this model is that the out-of-plane and in-plane dynamics are decoupled from each other, and the out-of-plane subsystem involves only a second-order harmonic oscillator, making it easier to analyze analytically than the coupled in-plane subsystem. Seo and Wiesel exploited this structure in their study of low-thrust optimal control for satellite formations, where they considered minimum-time control for the full relative dynamics and minimum-fuel control for the out-of-plane harmonic motion~\citep{force1,force2}. For the in-plane problem considered here, radial and along-track motion remain coupled, so a fixed-time minimum-fuel transfer must satisfy four terminal boundary conditions through two bounded control components. The two switching functions are different components of the same costate vector, and their relation is central to the switch-count and duration results below.

\cite{carter1984fuel,carter1987singular} studies fuel-optimal maneuvers under a Euclidean thrust bound, where full thrust can change direction continuously. He shows that the magnitude of the primer vector crosses one at most six times per orbital period, and he classifies the singular trajectories. \cite{pontani2015} use this bound for finite-thrust maneuvers in Euler--Hill relative motion, and \cite{pontani2015symmetry} studies symmetry properties of the corresponding thrust/coast sequences. \cite{hajek1977roots,hajek1979} provides general root bounds and finite bang/off existence results for linear systems. The geometry of bounded moment sets developed by \citet[Ch.~VII]{kreinnudelman1977} provides the tools we use to replace a singular control by a bang/off control with the same terminal moments.

For numerical HCW transfers, \cite{gazzino2018} optimize switching times for prescribed firing sequences. Junction equations that are linear in the costates once the event times are given have been used for this problem \citep{zanon2006,zanon2009}. \cite{carterbrient1991,carterbrient1992} develop finite-thrust optimality conditions and numerical solutions. \cite{deiuliis2023aas} study direct and indirect pseudospectral methods for finite-thrust formation reconfiguration and how discretization resolves switching profiles.

These studies leave open the joint switching structure under independent radial and along-track bounds. \cite{vanderha1991finite} bounds the along-track crossings alone, and the Euler--Hill bounds of \cite{pontani2015} assume a steerable thrust vector. To our knowledge, the restrictions that the common costate places on the pair of radial and along-track switch counts, and finite representatives of singular optima under these bounds, have not been characterized.

This paper studies the switching structure of fixed-time minimum-fuel in-plane HCW maneuvers. The switching and singular-arc analyses use the Pontryagin minimum principle (PMP)~\citep{leemarkus1967} in its normal form, with the fuel multiplier scaled to one. The two main analytical results are as follows.
\begin{enumerate}
\item \textbf{Switching counts and duration classification.}
Within one orbital period, a regular optimal control has at most four radial and six along-track switches, and hence at most ten. These component bounds leave 35 candidate pairs of radial and along-track switch counts, of which 11 have positive necessary duration bounds, 1 cannot be optimal within one orbit, and the remaining 23 occur at arbitrarily short durations for suitably chosen targets (see Switching Counts from the Common HCW Costate).
\item \textbf{Finite representatives of singular controls.}
When the radial switching function of an optimum remains at a threshold, the singular-axis control admits an optimal bang/off representative with at most three radial switches, while the along-track control has at most two. When the along-track switching function remains at a threshold, there is an optimal representative with zero radial thrust and at most four along-track switches. Both constructions preserve terminal moments, and therefore the terminal state and fuel (see Singular Controls and Finite-Event Representatives).
\end{enumerate}
Together with an existence and normality argument, the switching and singular-control results show that every target reachable within one orbital period has a minimum-fuel optimum with at most four radial and six along-track switches. We also develop a numerical construction built on these finite profiles (Section~\ref{sec:gm43_computation}).

\begin{figure}[t]
\makebox[\textwidth][c]{\includegraphics{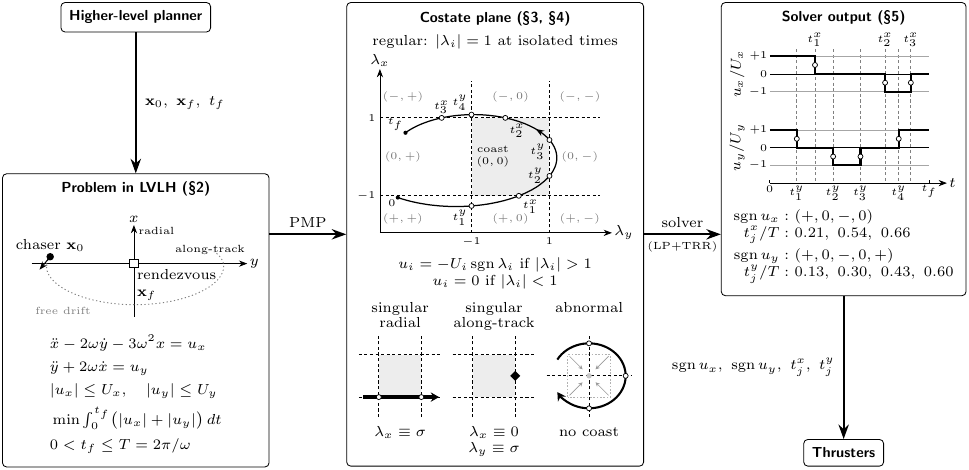}}
\caption{Overview of the paper through one regular example. A higher-level planner gives the duration ${t_f\approx0.75T}$, the initial state
${\mathbf x_0=[1.09,-13.19,-3.72,-3.07]^\top}$ and the final state
${\mathbf x_f=\mathbf 0}$ of an in-plane transfer in the LVLH frame, where
$T$ is the orbital period, in units with ${\omega=1}$ and ${U_x=U_y=1}$.
Here the transfer ends in rendezvous, and the dotted curve is the free drift
of the chaser. In the costate plane, the thresholds $\lambda_i=\pm1$ divide
the plane into nine cells labelled by the signs of $(u_x,u_y)$, and the
shaded cell is the coast region. The curve is the costate of this transfer,
and its $m=n_x+n_y=7$ threshold crossings are the control switches, labelled by
their switch times $t^x_j$ and $t^y_j$. The small planes below show the
other cases. In the two singular cases, a switching function stays on a
threshold. In the abnormal case, the fuel multiplier vanishes, so the coast
region shrinks to the origin and the controls switch where $\lambda_i$
changes sign. The right panel gives the optimal control of the same transfer
with the same switch times. The signs of $u_x$ and $u_y$ and the switch
times form the solver output that is sent to the thrusters. The solver is a
linear program (LP) followed by a trust-region reflective (TRR) correction.}
\label{fig:overview}
\end{figure}

Our numerical method starts from a discretized linear program (LP) that proposes control profiles and event times. We correct regular profiles by a trust-region reflective (TRR) iteration~\citep{colemanli1996}, using Gauss--Newton for profiles with fewer than four switches. We reconstruct singular-axis controls from their terminal moments and recover the costate from the switching conditions. A final whole-arc check verifies the PMP inequalities on every complete thrust and coast arc rather than at sampled times. For this convex problem, exact terminal feasibility and these conditions are sufficient for global minimum fuel.
Figure~\ref{fig:overview} shows the problem, the costate plane and the solver output for one transfer.

The scope of this paper is limited to local maneuver generation under the in-plane HCW dynamics, without collision avoidance, environmental perturbations and drag, thrust misalignment, or model uncertainty. The out-of-plane motion is decoupled and has been studied extensively by Seo and Wiesel~\citep{force1,force2}. The next section states the problem and its Pontryagin conditions.

\section{Problem Formulation}
\label{sec:problem}
\subsection{In-Plane Hill--Clohessy--Wiltshire Dynamics} \label{subsec:hcw_dynamics}

The motion of a spacecraft relative to a circular reference orbit is expressed in the LVLH frame, where \(x\) denotes the radial direction and \(y\) denotes the along-track direction. Under the HCW approximation, the cross-track motion is a forced harmonic oscillator decoupled from the in-plane motion~\citep{cw1960}, and the in-plane dynamics are
\begin{equation}
\ddot{x} - 2\omega \dot{y} - 3\omega^2 x = u_x,\qquad
\ddot{y} + 2\omega \dot{x} = u_y,
\label{eq:hcw_x}
\end{equation}
where \(\omega\) is the mean motion of the reference orbit and \(u_x\), \(u_y\) are the radial and along-track acceleration inputs. The Coriolis and gravity-gradient terms couple the two axes.

Define the in-plane state and control vectors as
\begin{equation}
\mathbf{x}(t) \coloneqq
\begin{bmatrix}
x(t) & y(t) & \dot{x}(t) & \dot{y}(t)
\end{bmatrix}^\top,
\qquad
\mathbf{u}(t) \coloneqq
\begin{bmatrix}
u_x(t) & u_y(t)
\end{bmatrix}^\top.
\end{equation}
Then the in-plane HCW equations in linear time-invariant (LTI) state-space form are
\begin{equation}
\dot{\mathbf{x}}(t) = A\mathbf{x}(t) + B\mathbf{u}(t),
\label{eq:inplane_state_space}
\end{equation}
with
\begin{equation}
A =
\begin{bmatrix}
0 & 0 & 1 & 0\\
0 & 0 & 0 & 1\\
3\omega^2 & 0 & 0 & 2\omega\\
0 & 0 & -2\omega & 0
\end{bmatrix},
\qquad
B =
\begin{bmatrix}
0 & 0\\
0 & 0\\
1 & 0\\
0 & 1
\end{bmatrix}.
\label{eq:inplane_matrices}
\end{equation}

\subsection{Minimum-Fuel Transfer Problem}
\label{subsec:minfuel_problem}

Consider a transfer over a prescribed maneuver time \(t_f>0\) between two in-plane relative states,
\begin{equation}
\mathbf{x}(0)=\mathbf{x}_0,
\qquad
\mathbf{x}(t_f)=\mathbf{x}_f .
\label{eq:boundary_conditions}
\end{equation}
The final time is fixed and is not treated as an optimization variable. The control inputs represent commanded accelerations along the LVLH \(x\) and \(y\) axes and are assumed to be bounded component-wise as
\begin{equation}
|u_x(t)| \leq U_{x}, \qquad |u_y(t)| \leq U_{y},
\qquad t\in[0,t_f],
\label{eq:control_bounds}
\end{equation}
where \(U_x>0\) and \(U_y>0\) are the maximum available accelerations along the two axes.

This box constraint models opposed thruster pairs that fire
independently along the radial and along-track LVLH axes, as in the
face-mounted thruster model of \cite{gazzino2018}. We
also assume equal effective exhaust velocity, nearly constant spacecraft mass, and no
coupling between the two firing channels. Propellant use is then
proportional to the total absolute acceleration effort:
\begin{equation}
J(\mathbf{u}) \coloneqq
\int_{0}^{t_f}
|u_x(t)|+|u_y(t)|\,dt .
\label{eq:fuel_cost}
\end{equation}
Thus, the in-plane fixed-time minimum-fuel transfer problem is stated as
\begin{align}
\min_{\mathbf{u}(\cdot)} \quad
& J(\mathbf{u})
\label{eq:optimal_control_problem}\\
\text{subject to} \quad
& \dot{\mathbf{x}}(t)=A\mathbf{x}(t)+B\mathbf{u}(t), \nonumber\\
& \mathbf{x}(0)=\mathbf{x}_0,\qquad \mathbf{x}(t_f)=\mathbf{x}_f, \nonumber\\
& |u_x(t)| \leq U_{x},\qquad |u_y(t)| \leq U_{y},
\qquad t\in[0,t_f]. \nonumber
\end{align}

The state transition matrix of system~\eqref{eq:inplane_state_space} is
$\Phi(t_2,t_1)=\exp[A(t_2-t_1)]$. Let $\mathbf k$ denote the
drift-subtracted terminal correction. Its components are the required
radial and along-track position and velocity changes. For prescribed
initial and terminal states, it satisfies
\begin{equation}
\mathbf k\coloneqq\mathbf x_f-\Phi(t_f,0)\mathbf x_0
=\int_0^{t_f}\Phi(t_f,t)B\mathbf u(t)\,dt.
\label{eq:drift_subtracted_terminal_map}
\end{equation}

\subsection{Pontryagin Conditions and Convex Sufficiency}
\label{subsec:pmp_switching}

To characterize the structure of the optimal control, we apply the PMP to the fixed-time problem~\eqref{eq:optimal_control_problem}. Let
\[
\boldsymbol{\lambda}(t)=
\begin{bmatrix}
\lambda_1(t) & \lambda_2(t) & \lambda_x(t) & \lambda_y(t)
\end{bmatrix}^{\top}
\]
denote the costate vector. The normal-form Hamiltonian is
\begin{equation}
H(\mathbf{x},\mathbf{u},\boldsymbol{\lambda})
\coloneqq
|u_x|+|u_y|
+
\boldsymbol{\lambda}^{\top}
\left(A\mathbf{x}+B\mathbf{u}\right).
\label{eq:hamiltonian}
\end{equation}
In this form of the PMP, the multiplier of the fuel objective
is nonzero and has been scaled to one. The regular and singular-arc
analyses below use this form. A minimum-fuel control exists for every
reachable target~\citep{leemarkus1967}. If its fuel multiplier vanished,
the control would push the terminal state as far as possible in one
direction, regardless of fuel. Formally, the PMP would make $\mathbf x_f$
minimize
$\boldsymbol\lambda(t_f)^\top\mathbf x$ over the convex
set of states reachable at $t_f$, so $\mathbf x_f$ would lie on
the boundary of that set~\citep{leemarkus1967}. Every target in its interior therefore has a
normal multiplier, and Proposition~\ref{prop:ten_switch} treats the
boundary.
The costate equations follow from the PMP relation
\(\dot{\boldsymbol{\lambda}}=-\partial H/\partial \mathbf{x}\).
Since the state-dependent part of the Hamiltonian is \(\boldsymbol{\lambda}^{\top}A\mathbf{x}\), this relation gives
\begin{equation}
\dot{\boldsymbol{\lambda}}(t)=-A^{\top}\boldsymbol{\lambda}(t).
\label{eq:costate_vector}
\end{equation}
Together with the state dynamics, the costate dynamics in Eq.~\eqref{eq:costate_vector} form a boundary-value problem. Hence, the initial costate variables are unknown parameters.

By the PMP pointwise minimization condition, the optimal control minimizes
the control-dependent part of the Hamiltonian,
\(\sum_{i\in\{x,y\}}\bigl(|u_i|+\lambda_i(t)u_i\bigr)\),
subject to~\eqref{eq:control_bounds} for almost every \(t\).
The Hamiltonian terms and control bounds separate by axis, so for \(i\in\{x,y\}\), the minimizing control satisfies
\begin{equation}
u_i^*(t)\in
\begin{cases}
\{-U_{i}\operatorname{sgn}\!\left(\lambda_i(t)\right)\},
& |\lambda_i(t)|>1,\\[2mm]
\{0\}, & |\lambda_i(t)|<1,\\[2mm]
[-U_{i},0], & \lambda_i(t)=1,\\[2mm]
[0,U_{i}], & \lambda_i(t)=-1.
\end{cases}
\label{eq:bang_off_bang_law}
\end{equation}
We call a control regular if neither switching function stays at
$\pm1$ on an open time interval. Each axis of a regular control is then
bang-off-bang, and $u_i$ switches exactly when $\lambda_i$ crosses $\pm1$.
If $\lambda_i$ stays at $\pm1$ on an interval, that axis has a singular arc,
on which condition~\eqref{eq:bang_off_bang_law} allows intermediate
thrust magnitudes. The section on singular controls treats this case.
This characterization does not by itself determine the number or timing
of the switches, as the initial costate is unknown.

\begin{proposition}[Sufficiency]\label{prop:sufficiency}
Let \(\mathbf u^*\) satisfy the dynamics, endpoint conditions, and
control bounds. If a solution \(\boldsymbol\lambda\) of
Eq.~\eqref{eq:costate_vector} makes \(\mathbf u^*\) minimize the
Hamiltonian almost everywhere, then \(\mathbf u^*\) has minimum fuel
among all admissible controls.
\end{proposition}
\begin{proof}
For any other feasible control \(\mathbf v\), pointwise minimization gives
\begin{equation}
J(\mathbf v)-J(\mathbf u^*)
\geq-\int_0^{t_f}\boldsymbol\lambda^\top B(\mathbf v-\mathbf u^*)\,dt
=-\left[\boldsymbol\lambda^\top(\mathbf x_v-\mathbf x_{u^*})\right]_0^{t_f}
=0.
\label{eq:convex_sufficiency_proof}
\end{equation}
The middle equality uses the state and adjoint equations, and the last
uses the shared endpoints.
\end{proof}

\subsection{Finite Switching Profiles and Terminal Conditions}
\label{subsec:gm45_profiles}
We now express the terminal conditions for a prescribed bang/off profile in terms of its switching times. This finite representation relates the required state transfer to the event times and thrust signs, while the PMP conditions specify the switching-function thresholds at each event.

For regular controls, condition~\eqref{eq:bang_off_bang_law} gives
bang/off segments. On each axis $i\in\{x,y\}$, $B_i$ is the corresponding column
of $B$ and $n_i$ counts the switches at times $0<t<t_f$. A simultaneous
change on both axes counts once in each. The $n_i+1$ segments are
indexed by $j=0,\ldots,n_i$.

Prescribe segment signs $s^i_0,\ldots,s^i_{n_i}\in\{-1,0,1\}$ and ordered event times
$0<t^i_1<\cdots<t^i_{n_i}<t_f$, with $t^i_0=0$ and
$t^i_{n_i+1}=t_f$. On segment $j$ the acceleration is
$u_i(t)=U_{i}s^i_j$ for $t\in(t^i_j,t^i_{j+1})$.
Equation~\eqref{eq:drift_subtracted_terminal_map}
then becomes the finite-arc terminal equation
\begin{equation}
\mathbf k=
\sum_iU_i\sum_{j=0}^{n_i}s^i_j
\int_{t^i_j}^{t^i_{j+1}}\Phi(t_f,t)B_i\,dt.
\label{eq:gm45_arc_terminal}
\end{equation}
All integrals in Eq.~\eqref{eq:gm45_arc_terminal} have closed HCW
forms. At event $j$ on axis $i$, define $\sigma^i_j\coloneqq s^i_{j-1}+s^i_j$. Exactly one adjacent segment is a
burn, so $\sigma^i_j\in\{-1,1\}$ and the PMP threshold condition is
$\lambda_i(t^i_j)=-\sigma^i_j$.

Although the initial costate is unknown, its common structure limits
$n_x$ and $n_y$, as we show next.
\section{Switching Counts from the Common HCW Costate}
\label{sec:gm43_bounds}

A regular control is bang-off-bang on each axis, and $u_i$ switches
exactly when $\lambda_i$ crosses $\pm1$. We consider a regular minimum-fuel control
with $0<t_f\leq T$, where $T=2\pi/\omega$ is one orbital period, and its
radial and along-track switch counts $n_x$ and $n_y$, with total
$m=n_x+n_y$. The two counts are not
independent. Both switching functions are components of the same costate,
and Eq.~\eqref{eq:gm43_phase_coupling} below shows that the slope of
$\lambda_y$ is set by the value of $\lambda_x$.

In this section, Proposition~\ref{prop:component_bounds} bounds each axis separately by
$n_x\leq4$ and $n_y\leq6$, which leaves 35 candidate pairs $(n_x,n_y)$.
The coupling then sorts these pairs, as summarized in
Table~\ref{tab:gm43_duration_allocations}.
Proposition~\ref{prop:durations} gives positive necessary durations for
11 pairs, Proposition~\ref{prop:exclusion} shows that $(3,6)$ cannot
occur in a regular optimum, and Proposition~\ref{prop:short} shows that
the remaining 23 pairs occur at arbitrarily short durations for suitably
chosen targets.

Define $[c_1\ c_2\ c_3\ c_4]^\top\coloneqq\boldsymbol\lambda(0)$. Solving
Eq.~\eqref{eq:costate_vector} gives the switching functions
\begin{equation}
\begin{aligned}
\lambda_x(t)&=-\frac{c_1}{\omega}\sin\omega t
+\frac{2c_2}{\omega}(\cos\omega t-1)
+c_3\cos\omega t+2c_4\sin\omega t,\\
\lambda_y(t)&=\frac{2c_1}{\omega}(1-\cos\omega t)
+c_2\left(3t-\frac{4}{\omega}\sin\omega t\right)
-2c_3\sin\omega t+c_4(4\cos\omega t-3).
\end{aligned}
\label{eq:gm45_switching_rows}
\end{equation}
For the switch-count analysis, we write both switching functions more compactly as
\begin{equation}
\lambda_x=a\sin\omega t+b\cos\omega t-2d,\qquad
\lambda_y=2a\cos\omega t-2b\sin\omega t+3d\omega t+e,
\label{eq:gm43_switch_functions}
\end{equation}
where
\begin{equation}
a\coloneqq2c_4-\frac{c_1}{\omega},\qquad
b\coloneqq c_3+\frac{2c_2}{\omega},\qquad
d\coloneqq\frac{c_2}{\omega},\qquad
e\coloneqq\frac{2c_1}{\omega}-3c_4.
\label{eq:gm43_switch_constants}
\end{equation}
The constants $a$, $b$, and $d$ appear in both functions, while $e$
appears only in $\lambda_y$. Equation~\eqref{eq:gm43_switch_constants} can be
solved uniquely for $c_1$ to $c_4$, so every choice of $a$, $b$, $d$, and $e$
gives an initial costate.

\begin{remark}
Each switch of a regular control is an isolated sign-changing crossing of one of the thresholds $\lambda_i=\pm1$ and the bounds below count these crossings.
\end{remark}

\subsection{At Most Four Radial and Six Along-Track Switches}
We first bound the switches on each axis separately, using only the form
of Eq.~\eqref{eq:gm43_switch_functions}.
\begin{proposition}\label{prop:component_bounds}
For a regular minimum-fuel control with $0<t_f\leq T$,
\begin{equation}
n_x\leq4,\qquad n_y\leq6,\qquad m\leq10.
\label{eq:gm43_event_bounds}
\end{equation}
\end{proposition}
\begin{proof}
The component bounds follow directly from Eq.~\eqref{eq:gm43_switch_functions}.
On $(0,T)$, a nonconstant radial sinusoid has at most two roots at each
threshold, while a constant $\lambda_x$ has no regular switches. The
derivative of $\lambda_y$ is a constant plus a sinusoid. If this derivative
is zero throughout the interval, $\lambda_y$ is constant and has no regular
switches. Otherwise, the derivative has at most two zeros, so Rolle's
theorem gives at most three roots of each along-track threshold,
which is the six-crossing bound of \cite{vanderha1991finite}.
Restricting to $(0,t_f)$ preserves these bounds.
\end{proof}

The component bounds leave $5\times7=35$ candidate count pairs,
including zero switches on either axis.
Table~\ref{tab:gm43_duration_allocations} records their duration
restrictions, which can be applied before the target is known. For example,
when $t_f\leq2\arccos(3/4)/\omega$, a solver does not need to try $(3,5)$, $(4,5)$,
$(4,6)$, or $(3,6)$.

\begin{table}[!htbp]
\centering
\caption{Duration restrictions for regular switch counts
over $0<t_f\leq T$. A dash means that the count pair occurs at arbitrarily
short positive durations for suitably chosen targets, and $\varnothing$
marks a pair excluded for regular optima. Roman numerals refer to the
necessary lower bounds of Proposition~\ref{prop:durations},
$t_f>0.115T$ (I), $t_f>0.159T$ (II), and $t_f>0.230T$ (III).}
\label{tab:gm43_duration_allocations}
\begin{tabular}{c@{\qquad}ccccccc}
\toprule
 & \multicolumn{7}{c}{Along-track count $n_y$}\\
Radial count $n_x$ & 0 & 1 & 2 & 3 & 4 & 5 & 6\\
\midrule
0 & -- & -- & -- & -- & -- & -- & --\\
1 & -- & -- & -- & -- & -- & II & II\\
2 & -- & -- & -- & -- & -- & II & II\\
3 & -- & -- & -- & I & I & III & $\varnothing$\\
4 & -- & -- & -- & I & I & III & III\\
\bottomrule
\end{tabular}
\end{table}

\subsection{Necessary Duration Bounds for 11 Count Pairs}
We now use the coupling between $\lambda_x$ and the slope of $\lambda_y$ to show
that 11 of the 35 pairs need a minimum transfer duration.
To express the three duration conditions and the short-transfer
constructions more simply, we use the orbital phase $\theta=\omega t$
and its final value $h=\omega t_f$. We denote differentiation with
respect to $\theta$ by a prime, so $\lambda_y'=\dot\lambda_y/\omega$. To expose
the shared sinusoid in Eq.~\eqref{eq:gm43_switch_functions}, we write
$a\sin\theta+b\cos\theta=A_\lambda\sin(\theta+\psi)$, where
$A_\lambda\geq0$ and $\psi$ are its amplitude and phase. The switching
functions and the derivative of $\lambda_y$ then become
\begin{equation}
\begin{aligned}
\lambda_x&=A_\lambda\sin(\theta+\psi)-2d,\\
\lambda_y&=2A_\lambda\cos(\theta+\psi)+3d\theta+e,\\
\lambda_y'&=-2A_\lambda\sin(\theta+\psi)+3d=-2\lambda_x-d.
\end{aligned}
\label{eq:gm43_phase_coupling}
\end{equation}

\begin{proposition}\label{prop:durations}
Let a regular minimum-fuel control have $0<t_f\leq T$.
\begin{enumerate}
\item[I.] If $n_x\geq3$ and $n_y\geq3$, then $t_f>\arccos(3/4)/\omega$.
\item[II.] If $n_x\geq1$ and $n_y\geq5$, then $t_f>1/\omega$.
\item[III.] If $n_x\geq3$ and $n_y\geq5$, then $t_f>2\arccos(3/4)/\omega$.
\end{enumerate}
\end{proposition}
\begin{proof}
The three conditions use this coupling in different ways. Condition I
uses radial crossings to make $\lambda_y$
monotone over a short interval. Condition II uses along-track
crossings to exclude radial switches. Condition III combines the radial
amplitude restriction with the separation of the stationary
points of $\lambda_y$.

\medskip\underline{Three or more radial switches.}\enspace
Suppose first that $h\leq\arccos(3/4)$.
If $n_x\geq3$, both radial thresholds are crossed and an interior
radial extremum occurs. Crossing both thresholds requires the full
sinusoidal range to extend below $-1$ and above $+1$. Thus
\begin{equation}
A_\lambda>1+2|d|.
\label{eq:gm43_radial_amplitude}
\end{equation}
If the interior radial extremum is a maximum, then
$\sin(\theta+\psi)=1$ there. Every phase in the transfer interval lies
within $h$ of this maximum, so
$\sin(\theta+\psi)\geq\cos h\geq3/4$. Consequently,
\begin{equation}
\lambda_y'=-2A_\lambda\sin(\theta+\psi)+3d
\leq-\frac32A_\lambda+3|d|
<-\frac32.
\label{eq:gm43_short_monotone}
\end{equation}
If the interior radial extremum is a minimum, then
$\sin(\theta+\psi)\leq-\cos h\leq-3/4$ throughout the interval.
The same amplitude bound gives
$\lambda_y'\geq3A_\lambda/2-3|d|>3/2$. Thus $\lambda_y$ is
strictly monotone in either case and $n_y\leq2$, so $n_x\geq3$ and
$n_y\geq3$ require $t_f>\arccos(3/4)/\omega$, giving Condition I.

\medskip\underline{Five or more along-track switches.}\enspace
Now suppose $0<h\leq1$ and $n_y\geq5$.
At least one along-track threshold is crossed three times, so Rolle's
theorem gives two distinct interior zeros of $\lambda_y'$. Since five
crossings rule out $\lambda_y'\equiv0$, this constant-plus-sinusoid
function has at most two zeros on an interval of length $h\leq1$.
Hence $\lambda_y$ has three monotone
pieces. At least two of them must each cross both thresholds to give
five switches. The two stationary values therefore lie on opposite
sides of both thresholds, and their absolute difference exceeds two.

Let the stationary phases be $\theta_1<\theta_2$, with separation
$\theta_2-\theta_1=2\alpha\leq h$. Because $2\alpha<\pi$, these roots
of $\lambda_y'=0$ surround a maximum or minimum of the shared sinusoid.
Equation~\eqref{eq:gm43_phase_coupling} then gives
$|3d|=2A_\lambda\cos\alpha$. Around a sine maximum, the cosine term in
$\lambda_y$ changes by $-4A_\lambda\sin\alpha$, while its linear term
changes by $4A_\lambda\alpha\cos\alpha$. The signs reverse around a
sine minimum. Thus the absolute difference between the stationary
values is $4A_\lambda(\sin\alpha-\alpha\cos\alpha)$. Moreover,
$\sin\alpha-\alpha\cos\alpha=\int_0^\alpha s\sin s\,ds
\leq\alpha^3/3$. Thus
\begin{equation}
2<4A_\lambda(\sin\alpha-\alpha\cos\alpha)
\leq\frac{4A_\lambda\alpha^3}{3}
\leq\frac{A_\lambda h^3}{6}.
\label{eq:gm43_short_excursion}
\end{equation}
Consequently $A_\lambda>12/h^3$. The same stationary-point identity
gives $|2d|=4A_\lambda\cos\alpha/3$. Since
$\cos\alpha\geq1-\alpha^2/2$ and $2\alpha\leq h$,
\[
\frac43\cos\alpha-1\geq\frac{2-h^2}{6}>0.
\]
Combining these bounds yields
\begin{equation}
|2d|-A_\lambda
=A_\lambda\left(\frac43\cos\alpha-1\right)
>\frac{2(2-h^2)}{h^3}\geq2,
\qquad 0<h\leq1.
\label{eq:gm43_short_radial_exclusion}
\end{equation}
The last inequality uses $2-h^2\geq h^3$ for $0<h\leq1$.
Because $\lambda_x=A_\lambda\sin(\theta+\psi)-2d$, its entire range is
below $-2$ when $d>0$ and above $2$ when $d<0$. It cannot cross a
radial threshold, so $n_x=0$. Consequently $n_x\geq1$ and
$n_y\geq5$ require $h>1$, or $t_f>1/\omega$, which is Condition II.

\medskip\underline{Three radial and five along-track switches together.}\enspace
Suppose $0<h\leq2\arccos(3/4)$, $n_x\geq3$, and $n_y\geq5$.
Each radial threshold has at most two crossings in this interval, so
three radial switches require both thresholds to be crossed. The
amplitude argument leading to inequality~\eqref{eq:gm43_radial_amplitude}
therefore applies and gives
\[
\left|\frac{3d}{2A_\lambda}\right|<\frac34.
\]
Because $n_y\geq5$, Rolle's theorem again gives two distinct interior
stationary points of $\lambda_y$. By
Eq.~\eqref{eq:gm43_phase_coupling}, their phases solve
\[
\sin(\theta+\psi)=\frac{3d}{2A_\lambda}.
\]
For $q=3d/(2A_\lambda)$, the solutions of $\sin u=q$ are
$u=\arcsin q$ and $u=\pi-\arcsin q$ modulo $2\pi$.
Their consecutive gaps alternate between $2\arccos|q|$ and
$2\pi-2\arccos|q|$.
Thus any two distinct stationary points are separated by at least
$2\arccos|q|>2\arccos(3/4)$. They cannot both fit inside an interval
of length $h\leq2\arccos(3/4)$. Therefore
\begin{equation}
\omega t_f>2\arccos(3/4),\qquad
t_f>\frac{\arccos(3/4)}{\pi}T.
\label{eq:gm43_ten_duration}
\end{equation}
The bound in Eq.~\eqref{eq:gm43_ten_duration} is Condition III for $(3,5)$, $(4,5)$, and $(4,6)$.
Conditions I--III give the 11 positive entries in
Table~\ref{tab:gm43_duration_allocations}.
\end{proof}

\begin{corollary}\label{cor:seven_switches}
For a regular minimum-fuel control with $0<t_f\leq T$,
\[
m\geq7\quad\Longrightarrow\quad t_f>\frac{\arccos(3/4)}{\omega}.
\]
\end{corollary}
\begin{proof}
If $n_x\geq3$, the bounds of Proposition~\ref{prop:component_bounds}
and $m\geq7$ imply $n_y\geq3$, so Condition I applies. Otherwise
$1\leq n_x\leq2$ and $n_y\geq5$, so Condition II gives
$t_f>1/\omega>\arccos(3/4)/\omega$.
\end{proof}

\subsection{Exclusion of \texorpdfstring{$(3,6)$}{(3,6)} from Regular Optima}
For one pair, the coupling rules out regular optima at every duration up
to one orbital period.

\begin{proposition}\label{prop:exclusion}
No regular minimum-fuel control with $0<t_f\leq T$ has
$(n_x,n_y)=(3,6)$.
\end{proposition}
\begin{proof}
We show that six along-track threshold crossings cannot coexist with exactly three radial threshold crossings under the common HCW costate.

Suppose $0<h\leq2\pi$, $n_x=3$, and $n_y=6$. Six sign-changing crossings of $\lambda_y=\pm1$ force two distinct interior
zeros $0<\theta_1<\theta_2<h$ of $\lambda_y'$. The points $\theta_1$ and $\theta_2$ are its only zeros
in $(0,2\pi)$. Each of the three monotone pieces can cross each
threshold at most once, so all six crossings require both thresholds
on every piece. Reversing both switching signs if needed, we take the
pieces to be decreasing, increasing, and decreasing. Hence
\[
\lambda_y(0)>1,\quad \lambda_y(\theta_1)<-1,\quad
\lambda_y(\theta_2)>1,\quad \lambda_y(h)<-1.
\]
Equation~\eqref{eq:gm43_phase_coupling} also defines $\lambda_y$ up to
$\theta=2\pi$. Its derivative has no further zero before that phase,
so $\lambda_y(2\pi)\leq\lambda_y(h)<\lambda_y(0)$. Since
$\lambda_y(2\pi)-\lambda_y(0)=6\pi d$, we obtain $d<0$.

At $\theta_1$ and $\theta_2$,
Eq.~\eqref{eq:gm43_phase_coupling} gives $\lambda_x=-d/2$. The same
identity shows that the two outer pieces have
$\lambda_x>-d/2$, while the middle piece has $\lambda_x<-d/2$.
Three radial crossings require both thresholds and
the amplitude bound $A_\lambda>1+2|d|=1-2d$ from
inequality~\eqref{eq:gm43_radial_amplitude}. If $d\leq-2$, the outer pieces
stay above $+1$, so $\lambda_x$ starts and ends above both thresholds.
Each threshold then has an even number of sign-changing crossings,
contrary to $n_x=3$. Thus $-2<d<0$. The outer pieces now lie above
$0$, so both $-1$ radial crossings are in the middle piece. The
remaining $+1$ crossing lies on one outer piece, so the other has none.
Denote this latter piece by $I$. At its stationary end $\lambda_x=-d/2<1$,
and $\lambda_x$ does not cross $+1$ on $I$, so $-d/2\leq\lambda_x\leq1$ on $I$.
The sinusoid has a global maximum above $1$
and minimum below $-1$, so $\lambda_x$ has no extremum within this range
and is monotone on $I$.

The along-track curve decreases across both thresholds on $I$, so its
decrease exceeds $2$. Set $z=\lambda_x(\theta)$ and
$z_*=\max_{\theta\in I}\lambda_x(\theta)\leq1$. Then
$-\lambda_y'=2z+d\geq0$, $dz=\lambda_x'(\theta)\,d\theta$, and
$|\lambda_x'|=\sqrt{A_\lambda^2-(z+2d)^2}$.
For $-d/2\leq z\leq1$, the amplitude bound gives
\[
A_\lambda^2-(z+2d)^2
>(1-2d)^2-(z+2d)^2
=1-z^2-4d(1+z)>1-z^2.
\]
Because
$\lambda_x$ is monotone on $I$, the change of variable gives increasing
$z$-limits after reversing them when $z$ decreases. The integrand is
nonnegative, so extending its upper limit from $z_*$ to $1$ gives
\[
\begin{aligned}
\int_I-\lambda_y'(\theta)\,d\theta
&=\int_{-d/2}^{z_*}
\frac{2z+d}{\sqrt{A_\lambda^2-(z+2d)^2}}\,dz\\
&\leq\int_{-d/2}^{1}
\frac{2z+d}{\sqrt{A_\lambda^2-(z+2d)^2}}\,dz\\
&<\int_{-d/2}^{1}\frac{2z}{\sqrt{1-z^2}}\,dz
=2\sqrt{1-\frac{d^2}{4}}<2.
\end{aligned}
\]
The decrease of $\lambda_y$ across $I$ is therefore below $2$, contradicting the required decrease above $2$.
Therefore, the pair $(3,6)$ cannot be optimal in a regular minimum-fuel control with
$0<t_f\leq T$.
\end{proof}

\subsection{Arbitrarily Short Attainability of 23 Count Pairs}
The remaining 23 pairs do not need a minimum duration, as we show by
construction.
\label{subsec:gm43_short_allocations}

\begin{proposition}\label{prop:short}
For each of the 23 pairs marked with a dash in
Table~\ref{tab:gm43_duration_allocations} and every sufficiently small
$t_f>0$, some target has a regular minimum-fuel control with exactly this
count pair, and every minimum-fuel control for that target has the same
counts.
\end{proposition}
\begin{proof}
For each of the 23 pairs, we construct a regular minimum-fuel control and
take its endpoint as the target. During a short transfer, each switching
function is close to a low-degree polynomial, so we choose a polynomial
with the desired threshold crossings and match the Taylor terms of one
switching function to it.

The constructions fall into three groups. In the first, a quadratic
radial curve and a monotone along-track curve give up to four radial and two along-track
switches. In the second, a quadratic along-track curve gives the opposite
range of counts. In the third, a cubic along-track curve supplies the remaining five- and
six-switch cases with no radial switches.
In each group, one exact HCW switching component converges to the chosen
polynomial as $h\to0^+$, and we count the switches of the other component directly.

We first rescale the time by $\xi=2t/t_f-1$, so $-1\leq\xi\leq1$.
After shifting the time origin to the midpoint,
Eq.~\eqref{eq:gm43_switch_functions} becomes
\begin{equation}
\begin{aligned}
\lambda_x&=a\sin(h\xi/2)+b\cos(h\xi/2)-2d,\\
\lambda_y&=2a\cos(h\xi/2)-2b\sin(h\xi/2)+3dh\xi/2+e.
\end{aligned}
\label{eq:gm43_rescaled_switching}
\end{equation}

The first two construction groups use the following five quadratics:
\begin{equation}
\begin{aligned}
P_0(\xi)&=\xi^2+2,&
P_1(\xi)&=\xi^2+\xi+\tfrac12,&
P_2(\xi)&=2\xi^2,\\
P_3(\xi)&=4(\xi-\tfrac13)^2-2,&
P_4(\xi)&=4\xi^2-2.
\end{aligned}
\label{eq:gm43_count_polynomials}
\end{equation}
On $[-1,1]$, $P_j$ has exactly $j$ crossings of the two levels
$\pm1$ combined. All crossings are transverse and interior, and neither
endpoint lies on a threshold. For example, $P_3$ has a minimum of
$-2$, starts above $+1$ and ends between the two thresholds. It crosses
$+1$ once and $-1$ twice.

\medskip\underline{Up to four radial and two along-track switches.}\enspace
We write a selected quadratic as
$P(\xi)=p_2\xi^2+p_1\xi+p_0$, with $p_2>0$, and choose
\begin{equation}
a=\frac{2p_1}{h},\qquad
b=-\frac{8p_2}{h^2},\qquad
d=\frac{b-p_0}{2}.
\label{eq:gm43_radial_polynomial_limit}
\end{equation}
The sine and cosine expansions give
$\lambda_x=P+O(h^2)$, uniformly with its first derivative in $\xi$.
The relation between the two switching functions gives, uniformly in $\xi$,
\[
\frac{d\lambda_y}{d\xi}
=\frac h2(-2\lambda_x-d)
=\frac{2p_2}{h}+O(h)\xrightarrow[h\to0^+]{}+\infty.
\]
Thus, the along-track curve is
strictly increasing, with
$\lambda_y(1)-\lambda_y(-1)=4p_2/h+O(h)\to+\infty$ as $h\to0^+$.
The coefficient $e$ shifts this curve vertically without changing $\lambda_x$.
We choose $e$ separately for each along-track count. If
$\lambda_y(-1)=2$, the increasing curve has 0 threshold crossings. If
$\lambda_y(-1)=0$, it crosses $+1$ once. Alternatively, if
$\lambda_y(0)=0$, its left endpoint lies below $-1$ and its right endpoint
lies above $+1$ for all sufficiently small positive $h$, so it crosses
each threshold once. Selecting $P_0$ through $P_4$ therefore gives every
pair with $0\leq n_x\leq4$ and $0\leq n_y\leq2$.

\medskip\underline{Up to two radial and four along-track switches.}\enspace
Independently choose one of the same quadratics for the along-track curve,
writing it as $P(\xi)=p_2\xi^2+p_1\xi+p_0$, $p_2>0$. Set
\begin{equation}
a=-\frac{4p_2}{h^2},\qquad
d=\frac{2b+2p_1/h}{3},\qquad e=p_0-2a.
\label{eq:gm43_along_polynomial_limit}
\end{equation}
For a fixed real number $r$, substituting
$b=-2(3p_2r+2p_1)/h$ gives, uniformly with the first derivatives,
\[
\lambda_y=P+O(h^2),\qquad
\lambda_x=\frac{2p_2}{h}(r-\xi)+O(h).
\]
For $r=2$, the radial curve stays above both thresholds and has no
switches. For $r=0$, it decreases through both thresholds and has two.
To obtain exactly one radial switch, instead choose
\begin{equation}
b=\frac{a\sin(h/2)+4p_1/(3h)}{\cos(h/2)-4/3}.
\label{eq:gm43_one_radial_choice}
\end{equation}
The resulting radial curve satisfies $\lambda_x(-1)=0$. Its derivative is
$-2p_2/h+O(h)<0$ and its right endpoint tends to $-\infty$.
Therefore, only the $-1$ threshold is crossed. Since this choice still has
$b=O(1/h)$, the limit $\lambda_y\to P$ is unchanged.
These choices give every pair with $0\leq n_x\leq2$ and
$0\leq n_y\leq4$.

\medskip\underline{Five or six along-track switches without radial switches.}\enspace
For a cubic $Q(\xi)=q_3\xi^3+q_2\xi^2+q_1\xi+q_0$, $q_3>0$,
we set
\begin{equation}
a=-\frac{4q_2}{h^2},\qquad b=\frac{24q_3}{h^3},\qquad
d=\frac{2b+2q_1/h}{3},\qquad e=q_0-2a.
\label{eq:gm43_cubic_limit}
\end{equation}
Then $\lambda_y=Q+O(h^2)$ with its first derivative, while
$\lambda_x=-8q_3/h^3+O(1/h)$ stays below both thresholds.
For six along-track switches, we use $Q(\xi)=8\xi^3-6\xi$, whose endpoint
values are $-2$ and $2$ and whose interior extrema are $2$ and $-2$.
Each of its three monotonic pieces crosses both thresholds.
For five switches, we use the same cubic in the shifted argument
$\mu=(19\xi+1)/20$, that is, $Q(\xi)=8\mu^3-6\mu$. Here, the left endpoint is $-0.432$, the two interior extrema remain
$2$ and $-2$, and the right endpoint is $2$. The three monotone
pieces now have one, two and two crossings.

\medskip\underline{Preservation of the counts at positive duration.}\enspace
The limiting polynomials cross the thresholds transversely in the
interior. Uniform convergence of the switching functions and their
first derivatives preserves one crossing near each polynomial root.
Away from these roots, the polynomials remain separated from both
thresholds, so uniform convergence excludes further crossings. The
other component has the stated count by monotonicity or by separation
from both thresholds.

The first two groups each give 15 pairs and share nine. The third
adds two, giving $15+15-9+2=23$. For every sufficiently small positive
$h$, these exact common-costate curves define regular minimum-fuel
controls for their reached endpoints by
Proposition~\ref{prop:sufficiency}. Any other optimum for the same
target turns inequality~\eqref{eq:convex_sufficiency_proof} into an
equality, so it also minimizes the Hamiltonian with the same costate
almost everywhere. This minimizer is unique away from isolated
crossings, so every optimum for that target has the same switch counts.
\end{proof}
\section{Singular Controls and Finite-Event Representatives}
\label{sec:gm43_singular}

Propositions~\ref{prop:component_bounds}--\ref{prop:short} apply to
regular controls. If a switching function persists at a threshold,
condition~\eqref{eq:bang_off_bang_law} allows a range of minimizing
amplitudes on that axis, and many amplitudes can be optimal for the same
target. We show that such an optimum still has a bang/off representative
with the same target and fuel, at most three radial and two along-track switches
for radial singularity, and at most four along-track switches for
along-track singularity.

\subsection{Persistent Thresholds}

The sine, cosine, and linear-time terms in Eq.~\eqref{eq:gm43_switch_functions}
imply that equality to a threshold on any open interval persists
throughout the horizon. For $\sigma\in\{-1,1\}$, they give the
following singularity class on each axis:
\begin{equation}
\begin{aligned}
\lambda_x\equiv\sigma &: &a=b=0,\quad d=-\sigma/2,&\quad
\lambda_y=e-3\sigma\omega t/2,\\
\lambda_y\equiv\sigma &: &a=b=d=0,\quad e=\sigma,&\quad
\lambda_x\equiv0.
\end{aligned}
\label{eq:gm43_singular_classes}
\end{equation}
Substituting a persistent threshold into the control terms of
Eq.~\eqref{eq:hamiltonian} gives $H_i(u_i)=|u_i|+\sigma u_i$, which
attains its minimum zero for every $u_i$ between $0$ and
$-\sigma U_{i}$. Thus every minimizing control on this axis can
be written as
\begin{equation}
u_i(t)=-\sigma U_{i}v(t),\qquad 0\leq v(t)\leq1.
\label{eq:gm43_singular_amplitude}
\end{equation}
The threshold sign fixes the physical thrust direction, so all burns
on a persistent-threshold axis have the same sign, while the target
determines which amplitudes $v$ are feasible.

The classes in~\eqref{eq:gm43_singular_classes} impose a separate restriction
through the shared HCW costate. If $\lambda_x\equiv\sigma$, then
$\lambda_y$ is affine and crosses each threshold at
most once. Its control therefore has at most two switches. If
$\lambda_y\equiv\sigma$, then $\lambda_x\equiv0$, so the radial axis has no
firing. The two classes
cannot occur together because the first gives $\lambda_y$ a nonzero
slope.

For $\sigma=+1$ in the radial-singular class, the decreasing
$\lambda_y$ permits the along-track profiles $(-)$, $(0)$, $(+)$,
$(-,0)$, $(0,+)$, and $(-,0,+)$. For $\sigma=-1$, their signs
reverse. The two signs give 12 possible radial-singular profile pairs.
With $u_x=0$, the along-track-singular class gives two possible
profile pairs, one for each sign $\sigma$. The moment construction below
shows that either alternating radial profile $(0,\pm,0,\pm)$ or
$(\pm,0,\pm,0)$ can represent every admissible radial moment vector, with one
physical sign on its burns.

Replacing $v$ must preserve the terminal state and fuel, which depend on $v$ only
through a few of its integrals. In Eq.~\eqref{eq:drift_subtracted_terminal_map},
the radial column of $\Phi(t_f,t)B$ contains only $1$, $\sin\omega t$, and
$\cos\omega t$, while the along-track position row adds a term linear in $t$.
We therefore use the terminal and auxiliary moments
\begin{equation}
\begin{aligned}
M_0(v)&\coloneqq\int_0^{t_f}v(t)\,dt,
&M_1(v)&\coloneqq\int_0^{t_f}v(t)\sin\omega t\,dt,\\
M_2(v)&\coloneqq\int_0^{t_f}v(t)\cos\omega t\,dt,
&M_3(v)&\coloneqq\int_0^{t_f}t v(t)\,dt,
&M_4(v)&\coloneqq\int_0^{t_f}t^2v(t)\,dt.
\end{aligned}
\label{eq:gm43_singular_moments}
\end{equation}
The radial terminal contribution depends on $M_0(v),M_1(v),M_2(v)$,
and the along-track contribution additionally requires $M_3(v)$. The fuel on
either singular axis is $J_i=U_{i}M_0(v)$.
We use the auxiliary moment $M_3(v)$ for radial singularity and
$M_4(v)$ for along-track singularity.

\paragraph{Moment selection and root count.}
Let $V$ be the set of all amplitudes $v$ with $0\leq v(t)\leq1$ on $[0,t_f]$.
The moments in Eq.~\eqref{eq:gm43_singular_moments} form the extended
moment sets
\begin{equation}
\begin{aligned}
\mathcal A_x&=\{(M_0(v),M_1(v),M_2(v),M_3(v)):v\in V\},\\
\mathcal A_y&=\{(M_0(v),M_1(v),M_2(v),M_3(v),M_4(v)):v\in V\}.
\end{aligned}
\label{eq:gm43_singular_augmented_sets}
\end{equation}
These sets are compact and convex~\citep[Ch.~VII]{kreinnudelman1977}. We hold the moments needed for the
terminal state and fuel at their original values $M_j^\circ$ and select
one point from the
corresponding slice by
\begin{equation}
\begin{aligned}
\text{radial:}\quad&\min_{v\in V}M_3(v)
&&\text{subject to }M_j(v)=M_j^\circ,\quad j=0,1,2,\\
\text{along-track:}\quad&\min_{v\in V}M_4(v)
&&\text{subject to }M_j(v)=M_j^\circ,\quad j=0,1,2,3.
\end{aligned}
\label{eq:gm43_singular_auxiliary_selection}
\end{equation}
Since each slice is nonempty and compact, such a minimum exists.
Minimizing the last coordinate places the selected point of the
extended moment set on the boundary of $\mathcal A_i$, since an
interior point would permit a smaller last coordinate without changing
the earlier ones.
At a boundary point of a convex moment set, the supporting-hyperplane
property~\citep[Ch.~I]{kreinnudelman1977} gives nonzero coefficients
$\boldsymbol\kappa$ for which the selected minimizer $v_*$ maximizes a
weighted sum of moments over all $v\in V$:
\begin{equation}
\sum_j\kappa_jM_j(v)\leq\sum_j\kappa_jM_j(v_*)
\qquad (v\in V).
\label{eq:gm43_singular_support}
\end{equation}
The sum runs through $j=0$ to $3$ in the radial case and to $4$ in the
along-track case, so both sides have the form
$\int_0^{t_f}f_i(t)v(t)\,dt$, where
\begin{equation}
\begin{aligned}
f_x(t)&=\kappa_0+\kappa_1\sin\omega t+\kappa_2\cos\omega t+\kappa_3t,\\
f_y(t)&=\kappa_0+\kappa_1\sin\omega t+\kappa_2\cos\omega t
       +\kappa_3t+\kappa_4t^2.
\end{aligned}
\label{eq:gm43_auxiliary_functions}
\end{equation}
The coefficients $\boldsymbol\kappa$ select an amplitude
representative, while the physical PMP costate stays fixed. Since the functions
in each $f_i$ are linearly independent and $\boldsymbol\kappa\ne0$,
$f_i$ is a nonzero analytic function with only finitely many zeros
on the closed horizon.
Inequality~\eqref{eq:gm43_singular_support} imposes the optimality
condition that $v_*$ maximize $\int_0^{t_f}f_i(t)v(t)\,dt$ over $V$. The
amplitude bound $0\leq v\leq1$ therefore yields, almost everywhere,
\begin{equation}
v_*(t)=
\begin{cases}
1,&f_i(t)>0,\\
0,&f_i(t)\leq0,
\end{cases}
\qquad
u_i^*(t)=-\sigma U_{i}v_*(t).
\label{eq:gm43_singular_indicator}
\end{equation}
Because the zeros of $f_i$ have measure zero, we set $v_*=0$ there
without changing any moment, the terminal state, or fuel.

The selected control can change only at a zero of $f_i$. On
$0<t_f\leq T$, both $f_x'$ and $f_y''$ are a constant plus a sinusoid,
so each has at most two interior zeros when it is not identically zero.
One application of Rolle's theorem gives at most three zeros of $f_x$,
and two applications give at most four zeros of $f_y$. An identically
zero derivative gives a constant or affine $f_i$ with fewer zeros.
Combining these bounds with the
other-axis restrictions in~\eqref{eq:gm43_singular_classes}, the
upper-bound pairs for $(n_x,n_y)$ are $(3,2)$ in radial singularity and
$(0,4)$ in along-track singularity. Consequently, for
$0<t_f\leq T$, every persistent-threshold optimum with a normal multiplier has a bang/off optimal representative with
the same target and fuel and one of these switch bounds.

\paragraph{Alternating moment profiles.}
We now show that, for each thrust sign $\sigma$, either of two
complementary bang/off sequences reproduces every feasible singular-axis
moment vector. Define
$G_x(t)\coloneqq(1,\sin\omega t,\cos\omega t)$ and
$G_y(t)\coloneqq(1,t,\sin\omega t,\cos\omega t)$, and set
\begin{equation}
\mathcal M_i(t_f)\coloneqq\left\{\int_0^{t_f}G_i(t)v(t)\,dt:v\in V\right\}.
\label{eq:gm43_singular_moment_body}
\end{equation}
These moment bodies are compact and convex. They are also
full-dimensional, by the independence of the derivative columns shown
below.
In forward time, choose one sequence from each complementary pair:
\[
\begin{aligned}
v^x&=(0,1,0,1),&\quad \bar v^x&=(1,0,1,0),\\
v^y&=(0,1,0,1,0),&\quad \bar v^y&=(1,0,1,0,1).
\end{aligned}
\]
The first sequences have three radial and four along-track event times,
matching the numbers of independent moments in $G_x$ and $G_y$.
It suffices to show that $v^i$ covers the entire feasible moment set.
Its complement $\bar v^i$ then does too. For the singular axis $i$, let
$r_x=3$ and $r_y=4$ be these numbers of event times, and let
$\Delta_i\coloneqq\{0\leq t_1^i\leq\cdots\leq t_{r_i}^i\leq t_f\}$,
with $t_0^i=0$ and $t_{r_i+1}^i=t_f$. The moment vector generated by
these event times is
\begin{equation}
F_i(\mathbf t^i)\coloneqq
\sum_{j=0}^{r_i}v_j^i\int_{t_j^i}^{t_{j+1}^i}G_i(t)\,dt,
\qquad
\frac{\partial F_i}{\partial t_j^i}
=\bigl(v_{j-1}^i-v_j^i\bigr)G_i(t_j^i).
\label{eq:gm43_singular_moment_jacobian}
\end{equation}
The same root-count argument gives at most two distinct interior zeros
for a nonzero combination of the entries of $G_x$, and at most three
for $G_y$, when $t_f\leq T$. These entries therefore form Chebyshev
systems~\citep{karlinstudden1966}. Thus the derivative columns are
independent at distinct interior times, since otherwise one such
combination would vanish at all of these times, and the event times can
reach every nearby moment. When times merge or reach an endpoint, the
shorter profile switches only where a combination of the entries of
$G_i$ changes sign, so it gives a boundary moment. If an interior moment
were missing, a line from an attained interior moment to it would have
a first gap. Since $\Delta_i$ is compact, the moment at that gap is
still attained. Its event times could be neither distinct, where the
derivatives cover nearby moments, nor degenerate, by the preceding
observation. Thus $F_i(\Delta_i)=\mathcal M_i(t_f)$, which is the
principal-representation result of Markov and
Krein~\citep{karlinstudden1966,kreinnudelman1977}.

For $\bar v_j^i=1-v_j^i$ and
$G_i^{\rm tot}=\int_0^{t_f}G_i(t)\,dt$,
\begin{equation}
F_{\bar v^i}(\mathbf t^i)=G_i^{\rm tot}-F_i(\mathbf t^i).
\label{eq:gm43_singular_complement}
\end{equation}
If $0\leq v\leq1$, then also $0\leq1-v\leq1$, and the moments of $1-v$
are $G_i^{\rm tot}$ minus the moments of $v$. Thus $v\mapsto1-v$ sends
$\mathcal M_i$ onto itself, and for any target
moment vector $\boldsymbol\eta\in\mathcal M_i$,
$G_i^{\rm tot}-\boldsymbol\eta$ is feasible. Therefore, choose event
times at which $v^i$ gives $G_i^{\rm tot}-\boldsymbol\eta$.
Equation~\eqref{eq:gm43_singular_complement} then shows that
$\bar v^i$ gives $\boldsymbol\eta$ at those same times. Together with
the fixed-sign control in Eq.~\eqref{eq:gm43_singular_amplitude}, the
complement map gives a complementary pair of profiles for each of the two physical
thrust signs on the singular axis. The control on the other axis obeys the
restrictions in~\eqref{eq:gm43_singular_classes}.

The next result combines these switch bounds with the regular analysis and removes the restriction to a normal multiplier.

\subsection{Every Reachable Target}
\begin{proposition}\label{prop:ten_switch}
Let $0<t_f\leq T$, and let $\mathbf x_f$ be reachable from $\mathbf x_0$
at time $t_f$ under the control bounds. Then the transfer has a
minimum-fuel control with $n_x\leq4$ and $n_y\leq6$, and hence at most
ten switches.
\end{proposition}
\begin{proof}
A minimum-fuel control exists, and interior targets have a normal
multiplier, as noted after Eq.~\eqref{eq:hamiltonian}. For these targets,
Proposition~\ref{prop:component_bounds} and the singular representatives
above give the bound. For a boundary target, a supporting hyperplane of
the convex reachable set~\citep{leemarkus1967} gives a nonzero $\boldsymbol\lambda(t_f)$ for
which $\mathbf x_f$ minimizes $\boldsymbol\lambda(t_f)^\top\mathbf x$
over that set. Let
$\boldsymbol\lambda(t)=\Phi(t_f,t)^\top\boldsymbol\lambda(t_f)$. Every feasible control reaches $\mathbf x_f$, so it also
minimizes $\boldsymbol\lambda(t)^\top B\mathbf u$ almost everywhere and
is at full thrust wherever $\lambda_x$ or $\lambda_y$ is nonzero. These
functions have the form of Eq.~\eqref{eq:gm43_switch_functions} with
threshold zero. The along-track one cannot vanish identically, since that
forces $a=b=d=e=0$ and hence $\boldsymbol\lambda=\mathbf0$. If the radial
one does not vanish identically either, the feasible control is unique,
and the root counts in the proof of
Proposition~\ref{prop:component_bounds} give at most two radial and three
along-track switches. Otherwise $a=b=d=0$, the along-track control is at
full thrust with one sign throughout, and the radial control solves a
minimum-fuel problem of the same type with the radial terminal column
alone. The same argument then gives at most four radial switches.
\end{proof}

\section{Numerical Construction and Validation}
\label{sec:gm43_computation}

\begin{figure}[!tbp]
\centering
\begin{tikzpicture}[font=\footnotesize,>={Stealth[length=4pt]},
  box/.style={draw,rounded corners=2pt,align=center,inner sep=2.5pt,minimum height=9mm},
  arr/.style={->,semithick}]
\node[box,text width=22mm] (lp) at (0,0) {Linear program\\Eq.~\eqref{eq:gm43_lp_seed}};
\node[box,text width=20mm,right=4mm of lp] (cand) {Candidate\\profile};
\node[box,text width=36mm,anchor=west] (low) at ($(cand.east)+(7mm,0)$) {Regular, \mbox{$m<4$}\\Gauss--Newton, Eq.~\eqref{eq:gm45_arc_terminal}};
\node[box,text width=36mm,above=2mm of low] (reg) {Regular, \mbox{$m\geq4$}\\TRR, Eq.~\eqref{eq:gm45_full_event_system}};
\node[box,text width=36mm,below=2mm of low] (sing) {Persistent threshold\\damped Newton, Eq.~\eqref{eq:gm43_singular_moment_jacobian}};
\node[box,text width=25mm,anchor=west] (chk) at ($(cand.east-|low.east)+(6mm,0)$) {Terminal, junction,\\and whole-arc\\PMP checks};
\node[box,text width=16mm,right=6mm of chk] (acc) {Accepted\\within\\tolerances};
\draw[arr] (lp) -- (cand);
\draw[arr] (cand.east) -- ++(3mm,0) |- (reg.west);
\draw[arr] (cand.east) -- ++(3mm,0) |- (sing.west);
\draw[arr] (reg.east) -- ++(3mm,0) |- (chk.west);
\draw[arr] (sing.east) -- ++(3mm,0) |- (chk.west);
\draw[arr] (cand.east) -- (low.west);
\draw[arr] (low.east) -- (chk.west);
\draw[arr] (chk) -- node[above,font=\scriptsize]{pass} (acc);
\draw[arr] (reg.north) -- ++(0,3mm) -| node[pos=.25,above,font=\scriptsize]{endpoint arc vanishes, event removed} (cand.north);
\draw[arr] (chk.south) |- ($(sing.south)+(0,-4mm)$) -| (lp.south);
\node[below,font=\scriptsize] at ($(sing.south)+(0,-4mm)$) {rejected, next proposal or finer LP grid};
\end{tikzpicture}
\caption{Numerical construction of a minimum-fuel control.}
\label{fig:gm43_procedure}
\end{figure}
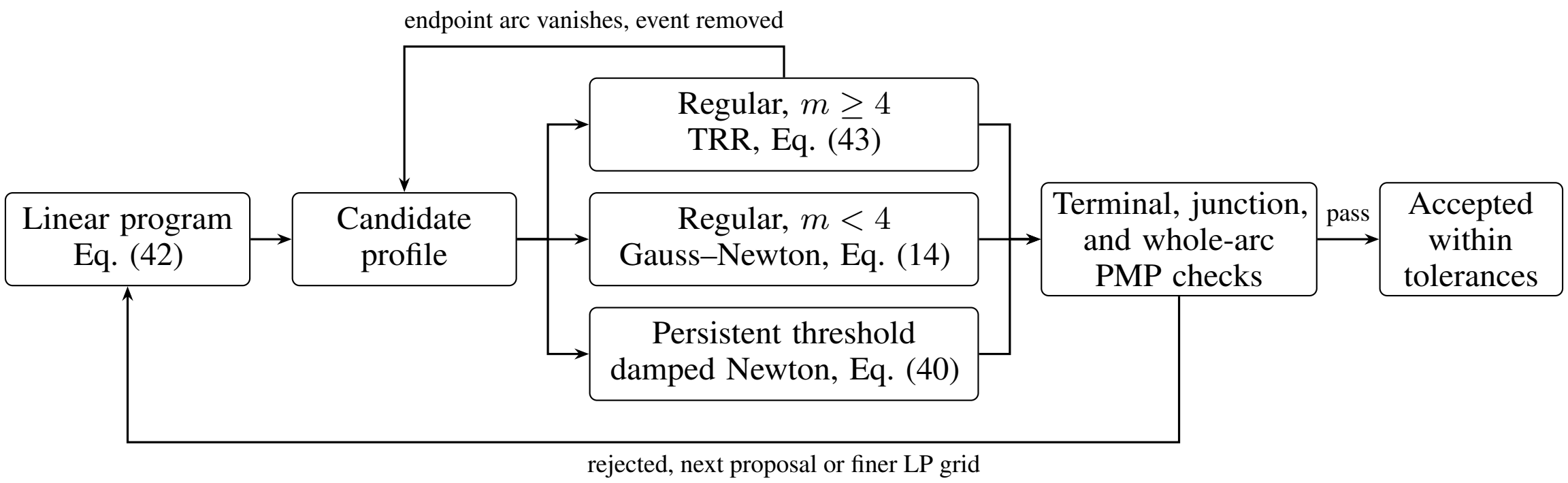

In this section, we compute minimum-fuel transfers within the finite
profile families above, following the procedure summarized in
Fig.~\ref{fig:gm43_procedure}.

\subsection{Linear-Program Initialization and Candidate Profiles}

We discretize the linear terminal map and absolute-value fuel to
initialize candidate thrust signs and event times. With
$K(t)\coloneqq\Phi(t_f,t)B$, divide $[0,t_f]$ into $N$ cells of length
$\delta\coloneqq t_f/N$ and set
$\bar K_j\coloneqq\int_{(j-1)\delta}^{j\delta}K(t)\,dt$.
For constant acceleration $\mathbf u_j=(u_{x,j},u_{y,j})^\top$ on cell
$j$, we solve
\begin{equation}
\min_{\mathbf u_1,\ldots,\mathbf u_N}
\delta\sum_{j=1}^N(|u_{x,j}|+|u_{y,j}|),\qquad
\sum_{j=1}^N\bar K_j\mathbf u_j=\mathbf k,\qquad
|u_{x,j}|\leq U_x,\quad|u_{y,j}|\leq U_{y}.
\label{eq:gm43_lp_seed}
\end{equation}
Exact HCW cell integrals and positive and negative control parts yield
a linear program. We first try a candidate profile whose arc sequence
and event times are obtained from the dual switching functions and
their crossings of $\pm1$. If this candidate is rejected, we try a
distinct profile obtained by grouping the primal control cells into
thrust and coast arcs. If no candidate is accepted, we solve the LP
again on a finer grid, from $N=128$ up to $N=4096$ cells.

\subsection{Continuous Correction and Validation}
\label{subsec:gm43_acceptance}

For a proposed regular profile, we collect its radial event times $t_j^x$
and along-track event times $t_j^y$, $m$ in total, in the vector $\mathbf t\in\mathbb R^m$. The terminal residual $\mathbf R_T(\mathbf t)$ is the
finite-arc contribution in Eq.~\eqref{eq:gm45_arc_terminal} minus the
prescribed correction $\mathbf k$. At each event $(i,j)$, the coefficients of $c_1,c_2,c_3,c_4$ in
Eq.~\eqref{eq:gm45_switching_rows} evaluated at $t_j^i$ form a row of
$M(\mathbf t)$, and the threshold $-\sigma_j^i$ enters $\mathbf b$. The
candidate equations are then
\begin{equation}
\begin{aligned}
\mathbf R_T(\mathbf t)
&\mathrel{:=} \sum_iU_i\sum_{j=0}^{n_i}s_j^i
\int_{t_j^i}^{t_{j+1}^i}\Phi(t_f,t)B_i\,dt-\mathbf k,\\
\mathbf R_T(\mathbf t)&=\mathbf0,
\qquad M(\mathbf t)\boldsymbol\lambda(0)=\mathbf b.
\end{aligned}
\label{eq:gm45_full_event_system}
\end{equation}
Some of the terminal equations can be solved in closed form for two or
three event times once the other times are fixed. We use these
identities to eliminate those times, taking the solution branch that
keeps all event times real, ordered, and inside $(0,t_f)$. TRR steps
then correct the remaining times within a bounded box. After each step
the eliminated times are recomputed, and the step is rejected if the
full set of event times is no longer ordered. When
$\operatorname{rank}M(\mathbf t)=4$, four independent junction rows
determine $\boldsymbol\lambda(0)$, and the remaining junction residuals enter the
TRR residual vector together with the reduced terminal equations. If
the correction fails while the first or last event time approaches
$0$ or $t_f$, we remove that event and solve and check the shorter
profile again.

Persistent-threshold candidates use the alternating moment profiles of
the preceding section. With the thrust sign fixed, three or four moment
equations determine the same number of event times. We solve them by
damped Newton continuation with the Jacobian in
Eq.~\eqref{eq:gm43_singular_moment_jacobian}, starting from the constant
amplitude $v=M_0^\circ/t_f$, which already has the required fuel. Regular
profiles with fewer than four switches are corrected by Gauss--Newton
iteration on the terminal equations alone, and their costate is then
obtained from the linear junction equations and the whole-arc
conditions.

We evaluate terminal states and fuel from constant-arc HCW formulas.
At every event the recovered costate must meet its threshold. On a burn
of sign $s$ we require $s\lambda_i\leq-1$, and on a coast
$|\lambda_i|\leq1$, throughout each complete arc. Evaluating the arc
endpoints and stationary points checks these inequalities without a
time grid. Together with exact endpoint feasibility, these conditions
imply global minimum fuel by Proposition~\ref{prop:sufficiency}. Because
they hold only to a numerical tolerance, we also record a relative
consistency bound, which adds the terminal residual weighted by the final
costate to the largest PMP violation integrated over both axes and the
whole transfer, and divides the sum by the fuel. The bound is zero when
the conditions of Proposition~\ref{prop:sufficiency} hold exactly.
For every returned control, the relative terminal residual, the junction
error, the whole-arc PMP violation, and the consistency bound stay below $10^{-8}$.

\section{Numerical Examples}
\label{sec:numerical}

We test the method on 55,296 feasible transfers with prescribed terminal
states and durations, and then on a sequence of waypoint maneuvers.
Throughout this section, $\omega=0.0011~\mathrm{rad\,s^{-1}}$ and
$U_x=U_y=5\times10^{-6}~\mathrm{m\,s^{-2}}$.

\subsection{Transfer Targets}

We use 96 durations from $T/1536$ to $T$,
of which 32 are logarithmically spaced below $T/20$ and 64 are
linearly spaced from $T/20$ to $T$. Position and velocity components are scaled by $U_{x}t_f^2$
and $U_{x}t_f$, respectively. At each duration, we draw 96
independent directions in these scaled coordinates by normalizing
four-component Gaussian vectors.
For each direction, we maximize target scale under the terminal map
and thrust bounds of LP~\eqref{eq:gm43_lp_seed}.
Scaling the resulting admissible control
gives feasible targets closer to the origin. Along each direction, fractions
$0.10$, $0.30$, $0.50$, $0.70$, $0.90$, and $0.99$ of the computed feasible ray limit
give 576 feasible targets per duration, or $96\times96\times6=55{,}296$ in total.

Every returned control is independently propagated to check its final
position and velocity, and no feasible target is excluded from the reported
statistics.

\subsection{Accuracy and Switching Structure}

Across the 55,296 transfers, the largest independently measured
relative terminal residual is below $5\times10^{-10}$, the largest
junction error and whole-arc PMP violation are below
$1.2\times10^{-9}$, and the largest relative consistency bound is below
$5.2\times10^{-9}$.

\begin{figure}[!bp]
\centering
\includegraphics[width=0.94\textwidth]{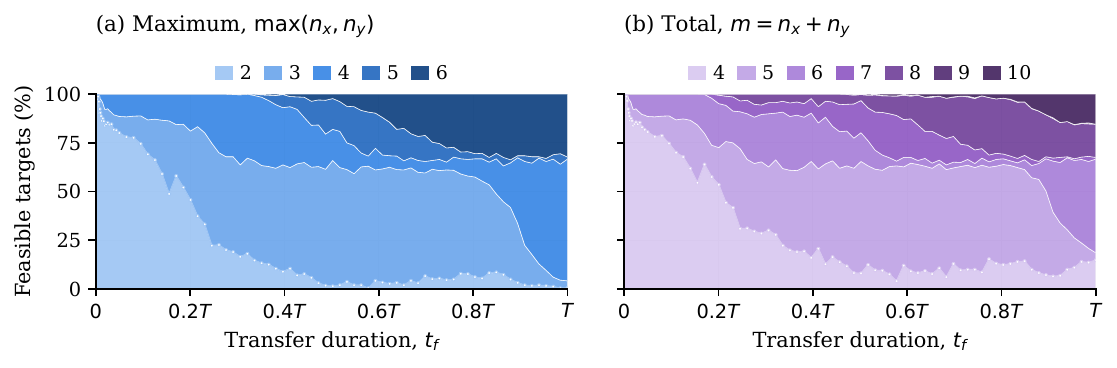}
\caption{Maximum component count (a) and total switch count (b) for the
same 55,296 feasible targets, with 576 at each duration.}
\label{fig:gm43_profiles}
\end{figure}

Figure~\ref{fig:gm43_profiles} groups the same controls by
$\max(n_x,n_y)$ in panel (a) and by the total $m$ in panel (b).
For example, $(0,6)$ and $(4,6)$ belong to the same maximum-six group
but have six and ten switches in total. The thin $m=9$ band in panel (b)
contains only 25 transfers, all with $(4,5)$, since $(3,6)$ cannot occur.
We observe 21 of the 35 count pairs. The most frequent are $(2,2)$
with 23,578 controls and $(2,3)$ with 13,512.
All transfers in count pairs restricted by
Table~\ref{tab:gm43_duration_allocations} have durations greater than
the corresponding bounds.

We also count a singular control through its computed bang/off
representative. Of the 55,296 controls, 257 use such representatives, of which
207 have $(n_x,n_y)=(0,4)$,
44 have $(3,1)$, and six have $(3,2)$. The 207 along-track-singular
representatives split almost evenly between the profiles $(0,+,0,+,0)$ and
$(0,-,0,-,0)$, and the 50 radial-singular ones use six of the 12
radial-singular profile pairs. Eight of the 14 sufficient profile-pair
families therefore appear in this sample.

\subsection{Individual and Sequential Waypoint Transfers}

\begin{figure}[!t]
\centering
\includegraphics[width=0.82\textwidth]{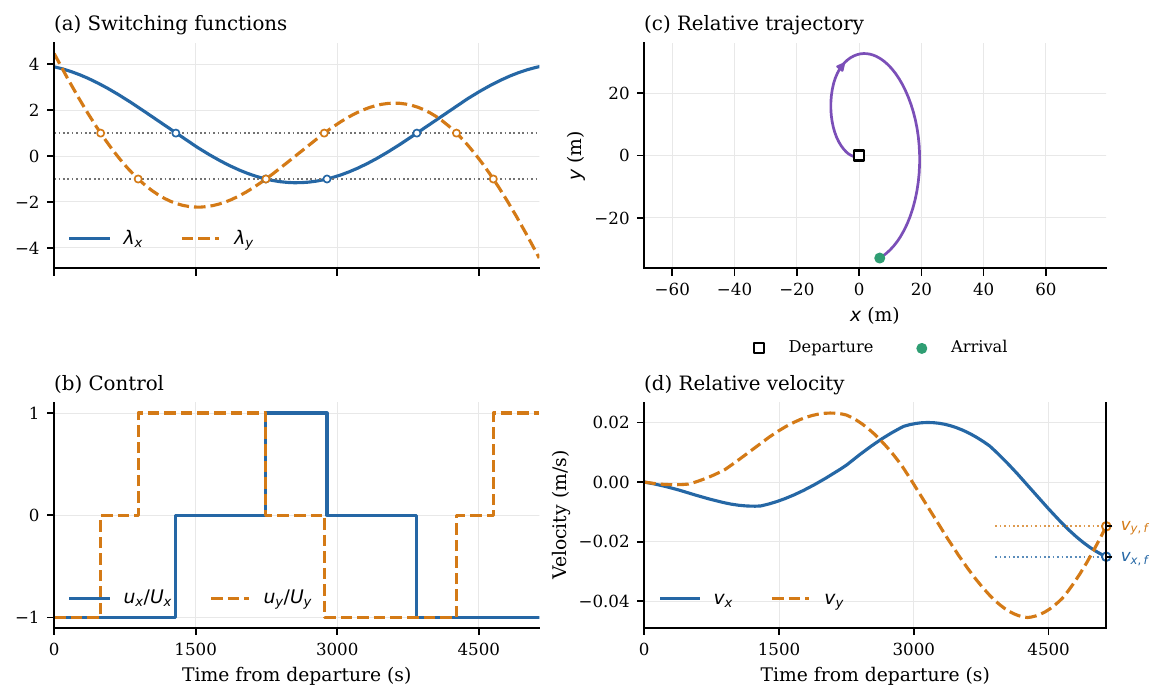}
\caption{First $0.9T$ transfer in Fig.~\ref{fig:gm43_waypoint}, with four
radial and six along-track switches. Open circles mark threshold
crossings. Panel (d) marks the target velocities.}
\label{fig:gm43_witness}
\end{figure}

\begin{figure}[!tbp]
\centering
\includegraphics[width=0.78\textwidth]{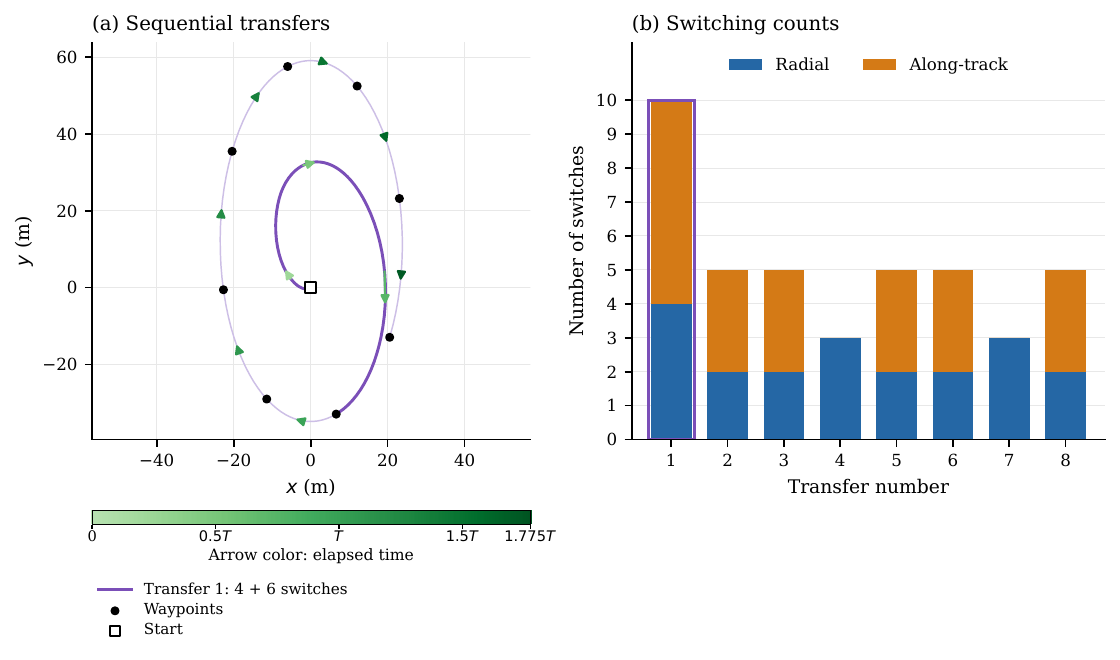}
\caption{Eight transfers. The first lasts $0.9T$ and the remaining
seven last $T/8$. Arrow shade shows elapsed time. The purple leg is
detailed in Fig.~\ref{fig:gm43_witness}. Stacked bars show component
switch counts.}
\label{fig:gm43_waypoint}
\end{figure}

We finally examine a sequence of eight waypoint transfers.
The first transfer lasts $0.9T$ and starts from rest at the origin.
The returned minimum-fuel control has four radial and six along-track
switches. Figure~\ref{fig:gm43_witness} shows its switching functions,
acceleration profiles, relative trajectory and velocity. The same
transfer is highlighted in Fig.~\ref{fig:gm43_waypoint}.

We generate the seven remaining targets by propagating the preceding
prescribed state for $T/8$ with constant acceleration
$0.12U_{x}(\cos[(j-1)\pi/3],\sin[(j-1)\pi/3])^\top$ for legs $j=2$
to $8$. These feasible reference controls specify the
waypoint states. Each leg starts from the arrival state obtained by
independently integrating the previous control and ends at the next
prescribed position and velocity. After the first leg, five legs have
$(2,3)$ switches and legs 4 and 7 have $(3,0)$, as
Fig.~\ref{fig:gm43_waypoint} shows for every leg.

All eight transfers also meet the terminal and whole-arc PMP
conditions. Within the HCW model, the maximum
position error is below $8\times10^{-12}~\mathrm m$ and the maximum
velocity error is below $7\times10^{-15}~\mathrm{m\,s^{-1}}$.
The largest relative consistency bound is below $2\times10^{-14}$.

\FloatBarrier
\section{Conclusions}
\label{sec:conclusion}

Every target reachable within one orbital period under the in-plane HCW
model has a minimum-fuel control with at most four radial and six
along-track switches. For regular optima, the common costate restricts
which radial and along-track counts can occur together. Because
Table~\ref{tab:gm43_duration_allocations} applies before the target is
known, it limits the profiles that a solver must try at a given
duration. For singular optima, moment selection gives bang/off
optima with at most five switches,
drawn from 12 radial-singular and two along-track-singular profile-pair
families.

All 55,296 tested transfers, including 257 singular representatives,
meet the endpoint and whole-arc PMP checks to within $10^{-8}$, and the
eight-leg waypoint sequence shows that the solver can be applied leg by
leg. The analytical bounds are established for horizons of at most one orbital period. The numerical acceptance test uses the normal form of the PMP, so the sample contains only targets strictly inside the reachable set; boundary targets, where the multiplier may be abnormal,
are covered by Proposition~\ref{prop:ten_switch} alone. Convergence of the numerical procedure is not proven, although it returned an accepted control for every sampled target. Extending the common-costate argument to other linear
relative-motion models, such as models for eccentric reference orbits, is
left for future work.

\section*{Funding Sources}

This work was supported in part by Horizon Europe MSCA COFUND Grant No. 101081645 and TÜBİTAK Grant No. 123C213 (both to B.B.).

\section*{Acknowledgements}
The authors would like to thank Utku Eren for fruitful discussions during the early stages of the work and for his constructive feedback on the manuscript. They also thank Abdurrahim Muratoglu for reviewing the draft.

\bibliographystyle{plainnat}
\bibliography{references_arxiv}

\begin{thebibliography}{28}
\providecommand{\natexlab}[1]{#1}
\providecommand{\url}[1]{\texttt{#1}}
\expandafter\ifx\csname urlstyle\endcsname\relax
  \providecommand{\doi}[1]{doi: #1}\else
  \providecommand{\doi}{doi: \begingroup \urlstyle{rm}\Url}\fi

\bibitem[Carter(1984)]{carter1984fuel}
Thomas~E. Carter.
\newblock Fuel-optimal maneuvers of a spacecraft relative to a point in circular orbit.
\newblock \emph{Journal of Guidance, Control, and Dynamics}, 7\penalty0 (6):\penalty0 710--716, 1984.
\newblock \doi{10.2514/3.19917}.

\bibitem[Carter(1987)]{carter1987singular}
Thomas~E. Carter.
\newblock Singular fuel-optimal space trajectories based on linearization about a point in circular orbit.
\newblock \emph{Journal of Optimization Theory and Applications}, 54\penalty0 (3):\penalty0 447--470, 1987.
\newblock \doi{10.1007/BF00940195}.

\bibitem[Carter and Brient(1991)]{carterbrient1991}
Thomas~E. Carter and John Brient.
\newblock Optimal bounded-thrust space trajectories based on linear equations.
\newblock \emph{Journal of Optimization Theory and Applications}, 70\penalty0 (2):\penalty0 299--317, 1991.
\newblock \doi{10.1007/BF00940628}.

\bibitem[Carter and Brient(1992)]{carterbrient1992}
Thomas~E. Carter and John Brient.
\newblock Fuel-optimal rendezvous for linearized equations of motion.
\newblock \emph{Journal of Guidance, Control, and Dynamics}, 15\penalty0 (6):\penalty0 1411--1416, 1992.
\newblock \doi{10.2514/3.11404}.

\bibitem[Clohessy and Wiltshire(1960)]{cw1960}
W.~H. Clohessy and R.~S. Wiltshire.
\newblock Terminal guidance system for satellite rendezvous.
\newblock \emph{Journal of the Aerospace Sciences}, 27\penalty0 (9):\penalty0 653--658, 1960.
\newblock \doi{10.2514/8.8704}.

\bibitem[Coleman and Li(1996)]{colemanli1996}
Thomas~F. Coleman and Yuying Li.
\newblock An interior trust region approach for nonlinear minimization subject to bounds.
\newblock \emph{SIAM Journal on Optimization}, 6\penalty0 (2):\penalty0 418--445, 1996.
\newblock \doi{10.1137/0806023}.

\bibitem[De~Iuliis et~al.(2023)De~Iuliis, Costigliola, and Casalino]{deiuliis2023aas}
Alessia De~Iuliis, Davide Costigliola, and Lorenzo Casalino.
\newblock Optimal control for formation reconfiguration problems using pseudospectral methods.
\newblock In \emph{33rd AAS/AIAA Space Flight Mechanics Meeting}, Austin, TX, 2023.
\newblock Paper AAS 23-260.

\bibitem[{Di Mauro} et~al.(2018){Di Mauro}, Bevilacqua, Spiller, Sullivan, and D'Amico]{DIMAURO2018311}
G.~{Di Mauro}, R.~Bevilacqua, D.~Spiller, J.~Sullivan, and S.~D'Amico.
\newblock Continuous maneuvers for spacecraft formation flying reconfiguration using relative orbit elements.
\newblock \emph{Acta Astronautica}, 153:\penalty0 311--326, 2018.
\newblock \doi{10.1016/j.actaastro.2018.01.043}.

\bibitem[Gazzino et~al.(2018)Gazzino, Arzelier, Cerri, Losa, Louembet, and Pittet]{gazzino2018}
Cl\'{e}ment Gazzino, Denis Arzelier, Luca Cerri, Damiana Losa, Christophe Louembet, and Christelle Pittet.
\newblock A minimum-fuel fixed-time low-thrust rendezvous solved with the switching systems theory.
\newblock \emph{Transactions of the Japan Society for Aeronautical and Space Sciences, Aerospace Technology Japan}, 16\penalty0 (7):\penalty0 628--634, 2018.
\newblock \doi{10.2322/tastj.16.628}.

\bibitem[H{\'a}jek(1977)]{hajek1977roots}
O.~H{\'a}jek.
\newblock On the number of roots of exp-trig polynomials.
\newblock \emph{Computing}, 18\penalty0 (2):\penalty0 177--183, 1977.
\newblock \doi{10.1007/BF02243627}.

\bibitem[H{\'a}jek(1979)]{hajek1979}
O.~H{\'a}jek.
\newblock $l_1$-optimization in linear systems with bounded controls.
\newblock \emph{Journal of Optimization Theory and Applications}, 29\penalty0 (3):\penalty0 409--436, 1979.
\newblock \doi{10.1007/BF00933143}.

\bibitem[Karlin and Studden(1966)]{karlinstudden1966}
S.~Karlin and W.~J. Studden.
\newblock \emph{Tchebycheff Systems: With Applications in Analysis and Statistics}, volume~15 of \emph{Pure and Applied Mathematics}.
\newblock Interscience Publishers, New York, NY, 1966.

\bibitem[Krein and Nudel'man(1977)]{kreinnudelman1977}
M.~G. Krein and A.~A. Nudel'man.
\newblock \emph{The Markov Moment Problem and Extremal Problems: Ideas and Problems of P. L. Chebyshev and A. A. Markov and Their Further Development}, volume~50 of \emph{Translations of Mathematical Monographs}.
\newblock American Mathematical Society, Providence, RI, 1977.
\newblock \doi{10.1090/mmono/050}.

\bibitem[Lee and Markus(1967)]{leemarkus1967}
E.~B. Lee and L.~Markus.
\newblock \emph{Foundations of Optimal Control Theory}.
\newblock Wiley, New York, NY, 1967.

\bibitem[Li et~al.(2023)Li, Zhang, Li, and He]{electronics12143034}
Xuehua Li, Lei Zhang, Zhijun Li, and Xingsuo He.
\newblock Application of the relative orbit in an on-orbit service mission.
\newblock \emph{Electronics}, 12\penalty0 (14):\penalty0 3034, 2023.
\newblock \doi{10.3390/electronics12143034}.

\bibitem[Luo et~al.(2014)Luo, Zhang, and Tang]{LUO20141}
Yazhong Luo, Jin Zhang, and Guojin Tang.
\newblock Survey of orbital dynamics and control of space rendezvous.
\newblock \emph{Chinese Journal of Aeronautics}, 27\penalty0 (1):\penalty0 1--11, 2014.
\newblock \doi{10.1016/j.cja.2013.07.042}.

\bibitem[Miller et~al.(2021)Miller, Walker, Agolli, and Dankanich]{miller2021survey}
Sara Miller, Mitchell~LR Walker, Jack Agolli, and John Dankanich.
\newblock Survey and performance evaluation of small-satellite propulsion technologies.
\newblock \emph{Journal of Spacecraft and Rockets}, 58\penalty0 (1):\penalty0 222--231, 2021.
\newblock \doi{10.2514/1.A34774}.

\bibitem[Opromolla et~al.(2018)Opromolla, Fasano, Rufino, and Grassi]{OPROMOLLA2018342}
Roberto Opromolla, Giancarmine Fasano, Giancarlo Rufino, and Michele Grassi.
\newblock Design of relative trajectories for in orbit proximity operations.
\newblock \emph{Acta Astronautica}, 145:\penalty0 342--356, 2018.
\newblock \doi{10.1016/j.actaastro.2018.01.062}.

\bibitem[Ortolano et~al.(2021)Ortolano, Geller, and Avery]{ortolano2021autonomous}
Nicholas Ortolano, David~K Geller, and Aaron Avery.
\newblock Autonomous optimal trajectory planning for orbital rendezvous, satellite inspection, and final approach based on convex optimization.
\newblock \emph{The Journal of the Astronautical Sciences}, 68\penalty0 (2):\penalty0 444--479, 2021.
\newblock \doi{10.1007/s40295-021-00260-5}.

\bibitem[Pontani(2015)]{pontani2015symmetry}
Mauro Pontani.
\newblock Symmetry properties of optimal relative orbit trajectories.
\newblock \emph{Mathematical Problems in Engineering}, 2015:\penalty0 286525, 2015.
\newblock \doi{10.1155/2015/286525}.

\bibitem[Pontani and Conway(2015)]{pontani2015}
Mauro Pontani and Bruce~A. Conway.
\newblock Minimum-fuel finite-thrust relative orbit maneuvers via indirect heuristic method.
\newblock \emph{Journal of Guidance, Control, and Dynamics}, 38\penalty0 (5):\penalty0 913--924, 2015.
\newblock \doi{10.2514/1.G000157}.

\bibitem[Seo(2004)]{force1}
J.~Seo.
\newblock \emph{Analytical Solution for Low-thrust Minimum Time Control of a Satellite Formation}.
\newblock PhD thesis, Air Force Institute of Technology, September 2004.

\bibitem[Seo and Wiesel(2004)]{force2}
John~S. Seo and William~E. Wiesel.
\newblock Low thrust control optimization for satellite formations: Analytical solution of the initial costates for time and fuel optimal harmonic motion.
\newblock In \emph{14th AAS/AIAA Space Flight Mechanics Conference}, volume 119, pages 2427--2445, Maui, HI, 2004.
\newblock Paper AAS 04-256.

\bibitem[Vallado(2022)]{vallado2022}
David~A Vallado.
\newblock \emph{Fundamentals of Astrodynamics and Applications}.
\newblock Microcosm Press, Torrance, CA, 5th edition, 2022.

\bibitem[van~der Ha(1991)]{vanderha1991finite}
Jozef~C. van~der Ha.
\newblock Analytical solutions for fuel-optimal rendezvous by finite-thrust arcs.
\newblock In \emph{42nd Congress of the International Astronautical Federation}, pages 1--20, Montreal, Canada, October 1991.
\newblock Paper IAF-91-363.

\bibitem[Yoo et~al.(2013)Yoo, Lee, Park, and Park]{YOO20131476}
Sung-Moon Yoo, Sangjin Lee, Chandeok Park, and Sang-Young Park.
\newblock Spacecraft fuel-optimal and balancing maneuvers for a class of formation reconfiguration problems.
\newblock \emph{Advances in Space Research}, 52\penalty0 (8):\penalty0 1476--1488, 2013.
\newblock \doi{10.1016/j.asr.2013.07.019}.

\bibitem[Zanon(2009)]{zanon2009}
Darren~J. Zanon.
\newblock \emph{Trajectory Planning for Spacecraft Formations Near Elliptical Orbits}.
\newblock PhD thesis, Cornell University, 2009.

\bibitem[Zanon and Campbell(2006)]{zanon2006}
Darren~J. Zanon and Mark~E. Campbell.
\newblock Optimal planner for spacecraft formations in elliptical orbits.
\newblock \emph{Journal of Guidance, Control, and Dynamics}, 29\penalty0 (1):\penalty0 161--171, 2006.
\newblock \doi{10.2514/1.7236}.

\end{thebibliography}

\end{document}